\documentclass[11pt]{article}

\usepackage{amssymb,amscd,array}
\usepackage[all]{xypic}

\let\ssection=\section
\renewcommand{\section}{\setcounter{equation}{0}\ssection}

\def\d{\delta}

\def\g{{\bf \mathrm g}}
\def\w{\mathbf{m}}
\def\om{\omega}
\def\r{\rho}
\def\a{\alpha}
\def\b{\beta}
\def\s{\sigma}
\def\vfi{\varphi}

\def\l{\lambda}
\def\m{\mu}
\def\n{\nabla}

\def\implies{\Rightarrow}

\newcommand{\am}{{\mathit a}}

\newcommand{\bbR}{\mathbb{R}}

\newcommand{\Diff}{\mathrm{Diff}}

\newcommand{\cF}{{\mathcal{F}}}

\newcommand{\Pol}{\mathrm{Pol}}
\newcommand{\cDp}{{\mathcal{D}}}

\newcommand{\SL}{\mathrm{SL}}

\newcommand{\Og}{\mathrm{O}}

\newcommand{\Vect}{\mathrm{Vect}}

\newcommand{\cqfd}{\hspace*{\fill}\rule{3mm}{3mm}}
\newcommand{\cqf}{\hspace*{\fill}\rule{2mm}{2mm}}
\begin{document}

\frenchspacing

\def\d{\delta}
\def\om{\omega}
\def\r{\rho}
\def\a{\alpha}
\def\b{\beta}
\def\s{\sigma}
\def\vfi{\varphi}
\def\l{\lambda}
\def\m{\mu}
\def\implies{\Rightarrow}

\oddsidemargin .1truein
\newtheorem{thm}{Theorem}[section]
\newtheorem{lem}[thm]{Lemma}
\newtheorem{cor}[thm]{Corollary}
\newtheorem{pro}[thm]{Proposition}
\newtheorem{ex}[thm]{Example}
\newtheorem{rmk}[thm]{Remark}
\newtheorem{defi}[thm]{Definition}
\title{Remarks on the Schwarzian Derivatives and the
Invariant Quantization by means of a Finsler Function}
\author{
Sofiane BOUARROUDJ\footnote{Supported by the Japan Society for the
Promotion of Science and by Action Concert\'ee de la Communaut\'ee
Fran\c{c}aise de Belgique.}\\}
\date{ {\footnotesize Department of Mathematics,}\\
{\footnotesize  Universit\'e Libre de Bruxelles,}\\
{\footnotesize  Campus Plaine CP 218, Bd. du Triomphe,}\\
{\footnotesize 1050-Brussels, Belgium.}} \maketitle
\begin{abstract} Let $(M,F)$ be a Finsler manifold. We construct a 1-cocycle
on $\Diff(M)$ with values in the space of differential operators
acting on sections of some bundles, by means of the Finsler
function $F.$ As an operator, it has several expressions: in terms
of the Chern, Berwald, Cartan or Hashiguchi connection, although
its cohomology class does not depend on them. This cocycle is
closely related to the conformal Schwarzian derivatives introduced
in our previous work. The second main result of this paper is to
discuss some properties of the conformally invariant quantization
map by means of a Sazaki (type) metric on the slit bundle
$TM\backslash 0$ induced by~$F.$
\end{abstract}
\section{Introduction}
The notion of {\it equivariant quantization} has been recently
introduced by Duval-Lecomte-Ovsienko in the papers
\cite{dlo,do,lo1}. The aim is to seek for an equivariant
isomorphism between the space of differential operators and the
corresponding space of symbols, intertwining the action of a Lie
group $G$ acting locally on a manifold $M$ -- see also \cite{bon,
b3, cmz, ga, lou} for related works. The computation was carried
out for the projective group $G=\SL (n+1,\bbR)$ in \cite{lo1}, and
for the conformal group $G=\Og(p+1,q+1),$ where $p+q=\mathrm
{dim}\,M$, in \cite{dlo,do}. It turns out that the
projectively/conformally equivariant quantization maps make sense
on any manifold, not necessarily flat, as shown in \cite{bor,
b3,do}. For instance, the conformally equivariant map has the
property that it does not depend on the rescaling of the (not
necessarily conformally flat) pseudo-Riemannian metric. The
existence of such maps induces naturally cohomology classes on the
group $\Diff(M)$ with values in the space of differential
operators acting on the space of tensor fields on $M$ of
appropriate types. These classes were given explicitly in
\cite{b1,b2, bo}, and interpreted as {\it projective} and {\it
conformal} multi-dimensional analogous to the famous {\it
Schwarzian derivative} (see \cite{b1,b2,bo} for more details).

A Riemannian metric is a particular case of more general functions
called {\it Finsler functions}. A Finsler function, whose
definition seems to go back to Riemann, is closely related to the
calculus of variation; it arises naturally in many context:
Physics, Mathematical ecology...(see\cite{aim}). The first main
result of this paper is to extend one of the two 1-cocycles
introduced by the author in \cite{b1} as conformal Schwarzian
derivatives, to the more general framework of Finsler structures.
The 1-cocycle can be built in terms of the Chern, Berwald, Cartan
or Hashiguchi connection. All of these connections are considered
as generalizations of the well-known Levi-Civita connection for
Riemannian structures, thereby the 1-cocycle in question coincide
with the conformal cocycle introduced in \cite{b1} when the
Finsler function $F$ is Riemannian. This 1-cocycle, thus, can be
exhibit in four ways accordingly to the used connection. This
property contrasts sharply with the case of projective structures
where the projectively invariant 1-cocycle has a unique expression
(cf. \cite{b2}).

The second part of this paper deals with the conformally invariant
quantization procedure. As the Finsler function gives rise to a
Riemannian metric, say $\w,$ on the slit bundle $TM \backslash 0,$
we shall apply the Duval-Ovsienko's quantization procedure through
the metric $\w.$ That means that we associate with functions on
the cotangent bundle of the manifold $TM \backslash 0$,
differential operators acting on the space of $\lambda$-densities
on $TM\backslash  0.$ The second main result of this paper is to
prove that, for almost all $\l$, this map cannot descend as an
operator acting on the space of $\l$-densities on $M$ even though
the Finsler function is Riemannian.
\section{Introduction to Finsler structures}
\label{fin}
We will follow verbatim the notation of \cite{bcs}. Let $M$ be a
manifold of dimension $n.$ A local system of coordinates $(x^i),
i=1,\ldots,n$ on $M$ gives rise to a local system of coordinates
$(x^i,y^i)$ on the tangent bundle $TM$ through\footnote{We will
use the convention of summation on repeated indices.}
$$
y=y^i\frac{\partial}{\partial x^i}, \quad \quad i=1,\ldots, n.
$$
\begin{defi}{\rm
A Finsler structure on $M$ is a function $F:TM\rightarrow
[0,\infty)$ satisfying the following conditions:

(i) the function $F$ is differentiable away from the origin;

(ii) the function $F$ is homogeneous of degree one in $y$, viz
$F(x,\l y)=\l F(x,y)$ for all $\l>0;$

(iii) the $n\times n$ matrix
$$\g_{ij}:=\frac{1}{2}\frac{\partial^2}{\partial y^i
\partial y^j}(F^2),$$
is positive-definite at every point of $TM\backslash 0.$}
\end{defi}
\begin{ex}{\rm
(i) Let $(M,\mathit{a})$ be a Riemannian manifold. The function
$F:=\sqrt{\mathit{a}_{ij}y^iy^j}$ satisfies the conditions (i),
(ii) and
(iii). In this case, the Finsler function is called {\it Riemannian}.\\
(ii)  Let $(M,\mathit{a})$ be a Riemannian manifold, and $\a$ be a closed
1-form on $M.$ We put $F:=\sqrt{\mathit{a}_{ij}y^iy^j}+\a_i y^i.$
One can prove that $F$ satisfies (i), (ii) and (iii) if and only if
$\|\a\cdot \a \|_{\mathit{a}}<1$ (see e.g. \cite{bcs}). In this case, $F$ is
Riemannian if and only if the 1-form $\a$ is identically zero.
}
\end{ex}
Denote by $\pi$ the natural projection $TM\backslash 0\rightarrow
M.$ The pull-back bundle of $T^*M$ with respect to $\pi$ is
denoted by $\pi^*(T^*M).$ The universal property implies that one
has a commutative diagram
\begin{equation}
\xymatrix{
T^*(TM\backslash 0)\ar[dr] \ar @{-->}[r]^u& \pi^{*}(T^*M) \ar[d]\ar[r]&T^*M\ar[d] \\
& TM\backslash 0 \ar[r]^\pi &M } \label{pi} \end{equation} The
components $\g_{ij}$ in (iii) of the definition above are actually
the components of a section of the pulled-back bundle $\pi^*
(T^{*}M)\otimes \pi^* (T^{*}M).$

The geometric object $\g$ in (iii) is called {\it fundamental}
tensor; it depends on $x$ and on $y$ as well. The fundamental
tensor is nothing but the Riemannian metric if $F$ is Riemannian.

The tensor
\begin{equation}
\label{car}
A:=A_{ijk}\,\, dx^i\otimes dx^j\otimes dx^k,
\end{equation}
where $A_{ijk}:=F/2\cdot  \partial \, \g_{ij}/\partial y^k,$
is called {\it Cartan} tensor.  It is symmetric on its three indices, and
defines a section of the pulled-back bundle $(\pi^* (T^{*}M))^{\otimes 3}.$
The Cartan tensor measures whether the Finsler function $F$ is Riemannian or
not.

The tensor
\begin{equation}
\label{hil}
\omega:=\omega_i\,\,dx^i,
\end{equation}
where $\omega_i:=\partial F/ \partial y^i,$ is called {\it Hilbert} form; it
defines a section of the pulled-back bundle $\pi^* (T^{*}M).$

Throughout this paper, indices are lowered or raised with respect
to the fundamental tensor $\g.$ For instance, the tensor whose
components are $A_{ij}^l$ stands for the tensor whose components
are $A_{ijk}\, \g^{kl}.$

We will also use the following notation: on the manifold
$TM\backslash 0,$ the index $i$ runs with respect to the basis $dx^i$ or
$\partial /\partial  x^i,$ and the index $\bar i$ runs with respect to the
basis $d y^i$ or $ \partial /\partial y^i.$
\section{The space of densities, the space of linear differential
operators and the space of symbols}
Let $M$ be an oriented manifold of dimension $n.$ Some backgrounds
are needed here to present our results. A thorough description of
all the forthcoming definitions can be found in \cite{dlo,do}.
\subsection{The space of densities and the space of linear differential
operators}
Let $(E,M)$ be a vector bundle over $M$ of rank $p.$ We define the
space of $\l$-densities of $(E,M)$ as the space of sections of the
line bundle $|\wedge^p E|^{\otimes \l}.$ Denote by $\cF_{\l}(M)$
the space of $\l$-densities associated with the bundle
$T^*M\rightarrow M$ and denote by $\cF_\l(\pi^* (T^* M))$ the
space of $\l$-densities associated with the bundle $\pi^*
(T^*M)\rightarrow TM\backslash 0$ (see (\ref{pi})). Both
$\cF_{\l}(M)$ and $\cF_\l(\pi^* (T^* M))$ are modules over the
group of diffeomorphisms $\Diff(M):$ for $f\in \Diff(M),$ $\phi\in
\cF_{\l}(M)$ and $\varphi \in \cF_\l(\pi^* (T^* M)),$ the actions
are given in local coordinates $(x,y)$ by
\begin{eqnarray}
\label{den}
f^*\phi&=&\phi\circ f^{-1}\cdot  ({J_{f^{-1}}})^{\lambda},\\
\label{den2} f^*\varphi&=&\varphi \circ \tilde f^{-1}\cdot
({J_{f^{-1}}})^{\lambda},
\end{eqnarray}
where $\tilde f$ is a lift of $f$ to $TM$ and $J_f=\left |Df/Dx \right |$ is
the Jacobian of $f$.

It is worth noticing that the forml{\ae} above do not depend on the choice of
the system of coordinates.

By differentiating these actions, one can obtain the actions of the Lie algebra
of vector fields $\Vect(M).$

Consider now ${\cal D}(\cF_\l(M),\cF_\mu (M)),$ the space of linear
differential operators
\begin{equation}
T:\cF_\l(M)\to\cF_\m(M).\nonumber
\label{Conv}
\end{equation}
The action of $\Diff(M)$ on ${\cal D}(\cF_\l(M),\cF_\mu (M))$ depends on the
two parameters $\l$ and $\m$; it is given by the equation
\begin{equation}
f_{\l,\m}(T)=f^*\circ T\circ {f^*}^{-1},
\label{Opaction}
\end{equation}
where $f^*$ is the action (\ref{den}) of $\Diff(M)$ on $\cF_\l(M)$.

Denote by ${\cal D}^2(\cF_\l(M),\cF_\mu (M))$ the space of second-order linear
differential operators with the $\Diff(M)$-module structure given by
(\ref{Opaction}). The space ${\cal D}^2(\cF_\l(M),\cF_\mu (M))$ is in fact a
$\Diff(M)$-submodule of ${\cal D}(\cF_\l(M),\cF_\mu (M)).$
\begin{ex}
{\rm The space of Sturm-Liouville operators
$\frac{d^2}{dx^2}+u(x): \cF_{-1/2}(S^1)\rightarrow \cF_{3/2}(S^1)$
on $S^1,$ where $u(x)\in\cF_{2}(S^1)$ is the potential, is a
submodule of $\cDp _{-\frac{1}{2},\frac{3}{2}}^2(S^1)$ (see
\cite{wi}). }
\end{ex}
Likewise, we define
${\cal D}(\cF_\l(\pi^* (T^* M)),\cF_\mu (\pi^* (T^* M))),$ the space of linear
differential operators
\begin{equation}
U :\cF_\l(\pi^* (T^* M))\to\cF_\m(\pi^* (T^* M)),\nonumber
\label{Conv2}
\end{equation}
with the action
\begin{equation}
f_{\l,\m}(U)=f^*\circ U \circ {f^*}^{-1}, \label{Opaction2}
\end{equation}
where $f^*$ is the action (\ref{den2}) of $\Diff(M)$ on
$\cF_\l(\pi^* (T^* M))$.

By differentiating the actions (\ref{Opaction}), (\ref{Opaction2}), one can
obtain the actions of the Lie algebra $\Vect(M).$

The formul{\ae} (\ref{Opaction}) and (\ref{Opaction2}) do not depend on the
choice of the system of coordinates.
\subsection{The space of symbols}
The space of symbols $\Pol (T^*M)$ is defined as the space of functions
on the cotangent bundle $T^*M$ that are polynomial on fibers.
This space is naturally isomorphic to the space
$\oplus_{p\geq 0}S\Gamma(TM^{\otimes p})$ of symmetric contravariant
tensor fields on $M.$

We define a one parameter family of $\Diff (M)$-module on the space of
symbols by
$$
\Pol_\d (T^*M):=\Pol (T^*M) \otimes \cF_{\d}(M).
$$
For $f\in \Diff(M)$ and $P\in \Pol_\d (T^*M),$ in local
coordinates $(x^i)$, the action is defined by
\begin{eqnarray}
\label{actsym}
f_{\d}(P)&=& f^*P\cdot (J_{f^{-1}})^{\d},
\end{eqnarray}
where $J_f=|Df/Dx|$ is the Jacobian of $f,$ and $f^*$ is the natural action
of $\Diff(M)$ on $\Pol (T^*M).$

We then have a graduation of $\Diff(M)$-modules given by
$$
\Pol_\d (T^*M)=\bigoplus_{k=0}^\infty \Pol_\d^k (T^*M),
$$
where $ \Pol_\d^k (T^*M)$ is the space of polynoms of degree $k$
endowed with the $\Diff(M)$-module structure (\ref{actsym}).
\begin{rmk}{\rm
As $\Diff(M)$-modules, the spaces $\Pol_\d (T^*M)$ and ${\cal
D}(\cF_\l(M),\cF_\mu(M))$ are not isomorphic (cf. \cite{dlo,do} ).
}
\end{rmk}
\section{Schwarzian derivative for Finsler structures}
Let $(M,F)$ be a Finsler manifold of dimension $n.$
\subsection{The Chern connection}
There exists a unique symmetric connection $D:
\Gamma(\pi^*(T^{*}M))\rightarrow \Gamma(\pi^*(T^{*}M)\otimes
T^*(TM\backslash 0))$ whose Christoffel symbols are given by
$$
\gamma^k_{ij}=\frac{1}{2}\g^{ks}\left (\frac{\partial \g_{si}}{\partial x^j}
+\frac{\partial \g_{sj}}{\partial x^i}-\frac{\partial \g_{ij}}{\partial x^s}
\right )
-\g^{ks}\left(\frac{N^m_j}{F}A_{msi}+\frac{N^m_i}{F}A_{msj}-
\frac{N^m_s}{F}A_{sij}\right ),
$$
where $A_{ijk}$ are the components of the Cartan tensor (\ref{car}),
$\g_{ij}$ are the components of the fundamental tensor and the components
$N_m^k$ are given by
\begin{equation}
\label{en}
N^k_m:=\frac{1}{4}\,\, \frac{\partial}{\partial y^m} \left (
\g^{ks}\left (\frac{\partial \g_{si}}{\partial x^j}
+\frac{\partial \g_{sj}}{\partial x^i}-\frac{\partial \g_{ij}}{\partial x^s}
\right )
\, y^j y^i\right )\cdot
\end{equation}
This connection is called {\it Chern connection} and has the following
properties (see \cite{bcs}):

(i) the connection 1-forms have no $dy$ dependence;

(ii) the connection $D$ is almost $\g$-compatible, in the sense
that $D_s(\g_{ij})=0$ and $D_{\bar s}(\g_{ij})=2\,A_{ijs};$

(iii) in general, the Chern connection is not a connection on $M;$
however, the Chern connection can descend to a connection on $M$
when $F$ is Riemannian. In that case, it coincides with the
Levi-Civita connection associated with the metric $\g.$
\subsection{A 1-cocycle as a Schwarzian derivative}
\label{pdef}
Since the connection 1-forms of the Chern connection have no $dy$
dependance, the difference between the two connections
\begin{equation}
\label{ell}
\ell(f):=f^*\gamma-\gamma,
\end{equation}
where $f\in \Diff(M),$ transforms under coordinates change as a section of
the bundle $\pi^* (T^{*}M)^{\otimes 2} \otimes \pi^* (TM).$
From the construction of the tensor (\ref{ell}), one can easily seen that
the map
$$
f\mapsto \ell(f^{-1})
$$
defines a non-trivial 1-cocycle on $\Diff(M)$ with values in
$S\Gamma (\pi^* (T^{*}M)^{\otimes 2}) \otimes \Gamma (\pi^* (TM)).$

Our main definition is the linear differential operator
${\cal A}(f)$ acting from
$S\Gamma (TM^{\otimes 2}) \otimes \cF_\d (M) $ to
$\Gamma (\pi^* (TM)) \otimes \cF_\d (\pi^* (T^*M)) $ defined by
\begin{equation}
\begin{array}{lcl}
{\cal A} (f)_{ij}^k&:=&\displaystyle
{f^*}^{-1} \left( \g^{sk}\,\g_{ij}\, D_s \right)
-\g^{sk}\,\g_{ij}\, D_s + (2-\delta n)\,\left(\ell(f)^k_{ij} -\frac{1}{n}
\, \mathrm{Sym}_{i,j}\,\delta_{i}^k \,\ell(f)^t_{tj}\right)\\[3mm]
&&\displaystyle +\g^{kl}\left ( \mathrm{Sym}_{i,j}\, \g_{sj}
B^s_{li}-\delta \, \g_{ij}\,B^t_{lt}\right ) -{f^{-1}}^{*} \left
(\g^{kl}\left ( \mathrm{Sym}_{i,j}\, \g_{sj} B^s_{li}-\delta\,
\g_{ij}\, B^t_{lt} \right )
    \right )\\[3mm]
&&\displaystyle
-(2-\delta n)\left ({f^{-1}}^{*} B^k_{ij}-B^k_{ij}-\frac{1}{n}
\mathrm{Sym}_{i,j}\, \delta^k_i ({f^{-1}}^{*} B^t_{jt}- B^t_{jt})\right ),
\label{MultiSchwar1}
\end{array}
\end{equation}
where we have put
$$
B^k_{ij}:=\left (A^{kr}_{i}\, \omega_j+A^{kr}_{j}\, \omega_i
-A^k_{ij}\, \omega^r-A^r_{ij}\, \omega^k -A^k_{is}\, A^{sr}_j
-A^k_{js}\,A^{sr}_i +A^{rk}_{u}\, A^{u}_{ij}\right ) d_r (\log F),
$$
to avoid clutter; $D$ is the Chern connection, $A_{ijk}$ are the
components of the Cartan tensor (\ref{car}), $\omega_i$ are the
components of the Hilbert form (\ref{hil}), $\ell(f)_{ij}^k$ are
the components of the tensor (\ref{ell}), $d_r:=u\circ d$ and $u$
is the map as in (\ref{pi}).
\begin{thm}
\label{main} {\rm (i)} The map $$f\mapsto {\cal A}(f^{-1}),$$
defines a 1-cocycle on $\Diff(M)$, non-trivial for all $\d \not=
2/n,$ with values in the space ${\cal D} (S\Gamma (TM^{\otimes
2})\otimes \cF_{\d}(M),
\Gamma(\pi^*(TM))  \otimes \cF_{\d}(\pi^*(T^*M)));$ \\
{\rm (ii)} The operator (\ref{MultiSchwar1}) does not depend on the rescaling
of the Finsler function $F$ by any non-zero positive function on $M;$ \\
{\rm (iii)} If $M:=\bbR^n$ and $F$ is Riemannian such that the metric
$\g$ is the flat metric, this operator vanishes on the conformal
group $\Og(n+1,1).$
\end{thm}
{\bf Proof.} Let us first explain how the contraction between the
tensor $D(P)$ and the tensor $\g^{-1}$ is permitted, for all $P\in
S\Gamma(TM^{\otimes 2})\otimes \cF_\d(M).$ Indeed, the tensor
$D(P)$ should take their values in $S\Gamma(TM^{\otimes 2})\otimes
\cF_{\d}(M)\otimes \Gamma (T^{*}(TM\backslash 0))$, from the
definition of the Chern connection. But taking into account that
the tensor $P$ lives in $S\Gamma(TM^{\otimes 2})\otimes
\cF_{\d}(M)$, the components $D_{\bar s}(P^{ij})=0,$ and the
components $D_{s}(P^{ij})$ behave under coordinates change as
components of a tensor in $\Gamma(TM^{\otimes 2})\otimes
\cF_{\d}(M)\otimes \Gamma (\pi^{*}(T^* M)).$ It follows therefore
that the contraction between $D_{s}(P^{ij})$ and $\g^{sk}$ makes a
sense.

To prove (i) we have to verify the 1-cocycle condition
$$
{\cal A}\left ((f\circ h)^{-1}\right )={f^*} {\cal A}(h^{-1})+
{\cal A}(f^{-1}), \quad \mbox{for all } f,h\in \Diff(M),
$$
where $f^*$ is the natural action on $\cDp (S\Gamma (TM^{\otimes
2}) \otimes \cF_{\d}(M), \Gamma(\pi^*(TM))  \otimes
\cF_{\d}(\pi^{*}(T^*M))).$ This condition holds because, in the
expression of the operator (\ref{MultiSchwar1}), $\ell$ is a
1-cocycle and the rest is a coboundary.

Let us proof that this 1-cocycle is not trivial for $\d\not=2/n$. Suppose
that there is a first-order differential operator
$A^k_{ij}=u^{sk}_{ij}D_s+v^k_{ij}$ such that
\begin{equation}
\label{cn} {\cal A}(f^{-1})=f^* A-A,
\end{equation}
It follows, by a direct computation, that
\begin{eqnarray}
\nonumber {f^*} v^k_{ij}-v^k_{ij}&=& (2-\delta
n)\left(\ell(f^{-1})^k_{ij} -\frac{1}{n}
\, \mathrm{Sym}_{i,j}\,\delta_{i}^k \,\ell(f^{-1})^t_{tj}\right)\\
\nonumber &&+\g^{kl}\left ( \mathrm{Sym}_{i,j}\, \g_{sj}
B^s_{li}-\delta \, \g_{ij}\,B^t_{lt}\right ) -{f}^{*} \left
(\g^{kl}\left ( \mathrm{Sym}_{i,j}\, \g_{sj} B^s_{li}-\delta\,
\g_{ij}\, B^t_{lt} \right )
    \right )\\
\nonumber &&-(2-\delta n)\left ({f}^{*}
B^k_{ij}-B^k_{ij}-\frac{1}{n} \mathrm{Sym}_{i,j}\, \delta^k_i
({f}^{*} B^t_{jt}- B^t_{jt})\right )\cdot
\end{eqnarray}
The right-hand side of this equation depends on the second jet of the
diffeomorphism $f,$ while the left-hand side depends on the first jet of
$f,$ which is absurd.

For $\d=2/n,$ one can easily see that the 1-cocycle (\ref{MultiSchwar1})
is a coboundary.

Let us prove (ii). Consider a Finsler function $\tilde F= \sqrt
\psi \cdot F,$ where $\psi$ is a non-zero positive function on
$M.$ Denote by $\tilde {\cal A}(f)$ the operator
(\ref{MultiSchwar1}) written by means of the function $\tilde F.$
To prove that $\tilde {\cal A}(f)={\cal A}(f)$ we proceed as
follows: we write down the tensors $\tilde D (P)$ and $\tilde
\ell(f)$ associated with the Finsler function $\tilde F$ in terms
of the tensor $D (P)$ and $\ell(f)$ associated with the Finsler
Function $F,$ then by replacing their expressions into the
explicit formula of the operator (\ref{MultiSchwar1}) we show that
the constants arising in the expression of the operator
(\ref{MultiSchwar1}) will annihilate the non-desired terms.

Let us first compare the Chern connections associated with the functions
$F$ and $\tilde F,$ namely
\begin{eqnarray}
\label{lien}
\tilde \gamma^k_{ij}&=&\gamma^k_{ij}+\frac{1}{2\psi}\left (\psi_i\, \d^k_j +
\psi_j \, \d^k_i
-\psi_t \,\g^{tk} \g_{ij}\right)\\
&&+\frac{1}{2\psi}\left (A^{kr}_{i}\, \omega_j
+A^{kr}_{j}\, \omega_i -A^k_{ij}\, \omega^r
-A^r_{ij}\, \omega^k  -A^k_{is}\, A^{sr}_j
-A^k_{js}\,A^{sr}_i +A^{rk}_{u}\, A^{u}_{ij}\right ) \psi_r,
\nonumber
\end{eqnarray}
where $\psi_r=\partial \psi/\partial x^r.$\\
From (\ref{lien}), a direct computation gives
\begin{eqnarray}
\label{bebe}
\tilde D_k P^{ij}&=&D_k P^{ij}+\frac{1}{2\psi}\left(  \mathrm{Sym}_{i,j}
P^{mi} \left(\psi_m \delta^j_k-\psi_t \,\g^{tj}\g_{km}\right) +(2-n\delta )\,
P^{ij}\psi_k\right),\\
\nonumber
&&+\mathrm{Sym}_{i,j} P^{sj} C^i_{ks}-\delta P^{ij} C^t_{tk} \\
\nonumber
\tilde \ell(f)^k_{ij}&=&\ell(f)^k_{ij}+{f^{-1}}^* \left (\frac{1}{2 \psi}
\left(\mathrm{Sym}_{i,j} \psi_i\delta^k_j -\psi_t\,
\g^{tk}\g_{ij}\right) \right )
-\frac{1}{2\psi}\left( \mathrm{Sym}_{i,j} \,\psi_i\, \delta^k_j -\psi_t\,
\g^{tk}\g_{ij}\right)\\
&&+ {f^{-1}}^* C^k_{ij}-C^k_{ij},
\nonumber
\end{eqnarray}
where we have put
$$
C^k_{ij}:=\frac{1}{2\psi}\left (A^{kr}_{i}\, \omega_j
+A^{kr}_{j}\, \omega_i -A^k_{ij}\, \omega^r -A^r_{ij}\, \omega^k
-A^k_{is}\, A^{sr}_j -A^k_{js}\,A^{sr}_i +A^{rk}_{u}\,
A^{u}_{ij}\right ) \psi_r,
$$
to avoid clutter; where $P^{ij}$ are the components of the tensor
$P\in S\Gamma(TM^{\otimes 2})\otimes \cF_{\d}(M).$\\
By substituting the formul{\ae} (\ref{bebe}) into (\ref{MultiSchwar1}) we get
by straightforward computation that ${\cal A} (f)=\tilde {\cal A} (f).$

Let us prove (iii). Suppose that $F$ is Riemannian, namely
$F=\sqrt{\g_{ij}y^iy^j}$. In that case, the Cartan tensor $A$ is
identically zero. The Chern connection $D$ can descend to a
connection on $M$ and coincide with Levi-Civita connection
associated with $\g.$ It follows that the operator
(\ref{MultiSchwar1}) turns into the form
\begin{equation}
\begin{array}{lcl}
{\cal A} (f)_{ij}^k&:=&\displaystyle
{f^*}^{-1} \left( \g^{sk}\,\g_{ij}\, D_s \right)
-\g^{sk}\,\g_{ij}\, D_s \\[2mm]
&&\displaystyle  + (2-\delta n)\,\left(\ell(f)^k_{ij} -\frac{1}{n}
\, \mathrm{Sym}_{i,j}\,\delta_{i}^k \,\ell(f)^t_{tj}\right)\cdot
\label{MultiSchwar1p}
\end{array}
\end{equation}
As the Chern connection and the fundamental tensor has no $y$
dependance, the operator (\ref{MultiSchwar1p}) will take its values in
$\Gamma(TM)\otimes \cF_{\d}(M),$ instead of
$\Gamma(\pi^{*}(TM))\otimes \cF_{\d}(\pi^*(T^*M)).$  The operator
(\ref{MultiSchwar1p}) is nothing but one of the conformally invariant
operators introduced in \cite{b3}. If, furthermore, $M$ is $\bbR^n$ and $\g$
is the flat metric then the operator (\ref{MultiSchwar1p}) vanishes on the
conformal group $\Og (n+1,1)$ (cf. \cite{b3}).\\
\cqfd
\begin{rmk}{\rm
(i) When the Finsler function $F$ is Riemannian, the formula
(\ref{MultiSchwar1p}) assures that the 1-cocycle $\cal A$ defined
here coincides with one of the conformally invariant operators
introduced in \cite{b1} as multi-dimensional conformal Schwarzian
derivatives. We refer to \cite{b1,b2,bo} for more explanations and
details concerning the relation between the classical Schwarzian
derivative and the projectively/conformally invariant operators
introduced in \cite{b1,b2,bo}.

(ii) The tensor whose components are $B^{kr}_{ij}$ coming out in
the formula (\ref{MultiSchwar1}) is only identically zero for
Finsler functions that are Riemannian. Indeed, if
$B^{kr}_{ij}\equiv 0$ then a contraction by the inverse of the
Hilbert form will give the equality $nA^k_{ij}\equiv 0.$ }
\end{rmk}

Now, how can we adjust the 1-cocycle $\cal A$ in order to take its
values in the space ${\cal D} (S\Gamma (TM^{\otimes 2})\otimes
\cF_{\d}(M), \Gamma(TM) \otimes \cF_{\d}(M))$, as for the
projectively/conformally invariant 1-cocycles of \cite{b1, b2,
bo}? A positive answer to this question can be given by demanding
an extra condition on the topology of $M.$ More precisely, suppose
that $M$ admits a non-zero vector fields - which is true when the
Euler Characteristic of $M$ is zero (cf. \cite{br}). One has
\begin{pro}
Let $\cal X$ be a fixed non-zero vector field on $M$ and denote by
$\breve{{\cal A}}(f)$ the operator obtained by substituting the
vector fields $\cal X$ into ${\cal A}(f)$ on the vertical
coordinates (namely $y$). The map
$$
f\mapsto \breve{{\cal A}}(f^{-1}),
$$
defines a 1-cocycle on $\Diff(M)$ with values in the space ${\cal
D} (S\Gamma (TM^{\otimes 2})\otimes \cF_{\d}(M), \Gamma(TM)
\otimes \cF_{\d}(T^*M)).$
\end{pro}
{\bf Proof.} Since the connection 1-forms of the Chern connection
do not depend on the direction of $dy,$ the evaluation by the
vector fields does not affect the 1-cocycle condition.  \cqf
\begin{rmk}{\rm
As a 1-cocycle, the cohomology class of $\breve{\cal A}$ does not
depend on the chosen vector fields. However, one has a family of
operators indexed by a family of non-vanishing vector fields on
$M.$ }
\end{rmk}
\subsection{The 1-cocycle ${\cal A}(f)$ in terms of the Berwald
connection}
There exists an other connection on the bundle
$\pi^*(TM)\rightarrow TM\backslash 0 $ called {\it Berwald}
connection. Its Christoffel symbols are given, in local
coordinates, by
\begin{equation}
{}^{\flat\!}\gamma^i_{jk}=\frac{\partial N^i_j}{\partial y^k},
\end{equation}
where the components $N^i_j$ are given as in (\ref{en}). Like the
Chern connection, this connection has no torsion (see \cite{mat}).

As in section \ref{pdef}, we define the following object
\begin{equation}
\label{psw}
{}^\flat\!\ell(f)=f^*{}^\flat\!\gamma-{}^\flat\!\gamma,
\end{equation}
where $f\in \Diff(M).$ As the connection 1-forms of the Berwald
connection have no $dy$ dependance, this object is actually a
section of the bundle $\pi^* (T^{*}M)^{\otimes 2} \otimes \pi^*
(TM).$

The 1-cocycle ${\cal A}(f)$ can be expressed in terms of the
Berwald connection as follows.
\begin{equation}
\begin{array}{lcl}
{}^\flat\!{\cal A} (f)_{ij}^k&:=&\displaystyle {f^*}^{-1} \left(
\g^{sk}\,\g_{ij}\, {}^\flat\! D_s \right) -\g^{sk}\,\g_{ij}\,
{}^\flat\! D_s + (2-\delta n)\,\left( {}^\flat\!\ell(f)^k_{ij}
-\frac{1}{n}
\, \mathrm{Sym}_{i,j}\,\delta_{i}^k \,{}^\flat\!\ell(f)^t_{tj}\right)\\[3mm]
&&\displaystyle +\g^{kl}\left ( \mathrm{Sym}_{i,j}\, \g_{sj}
{}^\flat\! B^s_{li}-\delta \, \g_{ij}\,{}^\flat\! B^t_{lt}\right )
-{f^{-1}}^{*} \left (\g^{kl}\left ( \mathrm{Sym}_{i,j}\, \g_{sj}
{}^\flat\! B^s_{li}-\delta \, \g_{ij}\,{}^\flat\! B^t_{lt} \right
)
    \right )\\[3mm]
&&\displaystyle -(2-\delta n)\left ({f^{-1}}^{*} {}^\flat\!
B^k_{ij}-{}^\flat\! B^k_{ij}-\frac{1}{n} \mathrm{Sym}_{i,j}\,
\delta^k_i ({f^{-1}}^{*} {}^\flat\! B^t_{jt}- {}^\flat\!
B^t_{jt})\right ), \label{pMultiSchwar1}
\end{array}
\end{equation}
where we have put
$$
{}^\flat\! B^k_{ij}:=\left ({}^\flat\! D_{\bar j} (A^{kr}_{i})\,
F+A^{kr}_{i}\, \omega_j +2\,A^{sk}_{j}\, \omega_i\right ) \, d_r
(\log F),
$$
to avoid clutter; ${}^\flat\! D$ is the covariant derivative
associated with the Berwald connection, $A_{ijk}$ are the
components of the Cartan tensor (\ref{car}), $\omega_i$ are the
components of the Hilbert form (\ref{hil}),
${}^\flat\!\ell(f)_{ij}^k$ are the components of the tensor
(\ref{psw}), $d_r:=u\circ d$ and $u$ is the map as in (\ref{pi}).

Theorem \ref{main} still holds for the operator ${}^\flat\! {\cal
A}(f).$ For the proof we proceed as in Theorem (\ref{main}). Part
(i) is obvious form the construction of the operator. Part (ii)
lead us to compare the Berwald connections associated with the
Finsler function $F$ and $\sqrt \psi\, F,$ respectively. The proof
then is a direct computation. Part (iii) results from the fact
that, as for the Chern connection, the Berwald connection
coincides with the Levi-Civita connection associated with a
Riemannian metric when $F$ is Riemannian.

It is worth noticing that, viewed as operators, the operator
${\cal A}(f)$ and ${}^\flat\!{\cal A}(f)$ are not equal; however,
they can be compared via the following definition.
\begin{defi} (see \cite{bcs}) {\rm
A Finsler manifold is called a {\it Landsberg} space if the tensor
whose components
$$
\dot{A}_{ijk}:=-\frac{1}{2}\, y_l\, \frac{\partial^2
N^l_i}{\partial y^j\partial y^k},
$$
where the components $N^l_i$ are as in (\ref{en}), is identically
zero. }
\end{defi}
\begin{pro}
\label{mk}{\rm If the Finsler manifold is a Landsberg space, the
operators
$${\cal A}(f)\equiv {}^\flat\!{\cal
A}(f),$$ for all $f\in \Diff(M).$ }
\end{pro}
{\bf Proof.} The proof results from the fact that the Berwald
connection and the Chern connection coincides when the Finsler
manifold is a Landsberg space (cf. \cite{bcs, mat}). \cqf

\begin{rmk}{\rm An obvious example of a Landsberg space is a
Riemannian manifold. More general examples of Finsler manifolds
that are Landsberg spaces but not Riemannian are provided in
\cite{bcs}. Proposition \ref{mk} shows that the operator ${\cal
A}(f)$ and ${}^\flat\! \!{\cal A}(f)$ does not coincide only for
Riemannian manifolds but also for some manifolds little more
general. }
\end{rmk}
\subsection{The 1-cocycle ${\cal A}(f)$ in terms of the Cartan
connection}
There exists an other connection on the bundle
$\pi^*(TM)\rightarrow TM\backslash 0 $ called {\it Cartan}
connection. Its Christoffel symbols are given, in local
coordinates, by
\begin{equation}
\label{mus}
\begin{array}{ll}
{}^{\natural\!}\gamma^i_{jk}=\displaystyle
\gamma^i_{jk}+A^i_{jt}\, \frac{N^t_k}{F},&\displaystyle
{}^{\natural\!}\gamma^i_{j\bar
k}=\frac{A^i_{jk}}{F},\\[3mm]
\displaystyle {}^{\natural\!}\gamma^i_{\bar j\bar k}=0, &
\displaystyle {}^{\natural\!}\gamma^i_{\bar j
k}=\frac{A^i_{jk}}{F},
\end{array}
\end{equation}
where $\gamma^i_{jk}$ are the Christoffel symbols of the Chern
connection, the components $N^i_j$ are as in (\ref{en}) and $
A^i_{jt}$ are defined in the section \ref{fin}. In
contradistinction with the Chern or Berwald connection, this
connection has the properties (cf. \cite{mat}):

(i) it has torsion;

(ii) the connection 1-forms do depend on the direction of $dy.$\\

As in section \ref{pdef}, we define the following geometrical
object
\begin{equation}
\label{ppsw} {}^\natural\!\ell(f)= {\tilde
f}^*{}^\natural\!\gamma-{}^\natural\! \gamma,
\end{equation}
where $\tilde f$ is a natural lift of $f\in \Diff(M).$ This object
takes its values in $\Gamma (\pi^{*}(TM)\otimes T^* (TM\backslash
0)\otimes T^* (TM\backslash 0)),$ in contrast with the previous
object $\ell(f)$ defined by means of the Chern or Berwald
connection. One takes the image of ${}^\natural\!\ell(f)$ by the
map $\mathrm{Id}\otimes u\otimes u,$ where $\mathrm{Id}$ is the
identity map and $u$ is the map defined in the diagram (\ref{pi}).
Let us still denote this tensor by ${}^\natural\!\ell(f).$

The 1-cocycle ${\cal A}(f)$ can be expressed in terms of the
Cartan connection as follows.
\begin{equation}
\begin{array}{lcl}
{}^\natural\!{\cal A} (f)_{ij}^k&:=&\displaystyle {f^*}^{-1}
\left( \g^{sk}\,\g_{ij}\,\,u\!\circ\! {}^\natural\! D_s \right)
-\g^{sk}\,\g_{ij}\,\, u\!\circ \!{}^\natural\! D_s + (2-\delta
n)\,\left( {}^\natural\!\ell(f)^k_{ij} -\frac{1}{n}
\, \mathrm{Sym}_{i,j}\,\delta_{i}^k \,{}^\natural\!\ell(f)^t_{tj}\right)\\[3mm]
&&\displaystyle +\g^{kl}\left ( \mathrm{Sym}_{i,j}\, \g_{sj}
{}^\natural\! B^s_{il}-\delta \, \g_{ij}\,{}^\natural\!
B^t_{tl}\right ) -{f^{-1}}^{*} \left (\g^{kl}\left (
\mathrm{Sym}_{i,j}\, \g_{sj} {}^\natural\! B^s_{il}-\delta \,
\g_{ij}\,{}^\natural\! B^t_{tl} \right )
    \right )\\[3mm]
&&\displaystyle -(2-\delta n)\left ({f^{-1}}^{*} {}^\natural\!
B^k_{ij}-{}^\natural\! B^k_{ij}-\frac{1}{n} \mathrm{Sym}_{i,j}\,
\delta^k_i ({f^{-1}}^{*} {}^\natural\! B^t_{tj}- {}^\natural\!
B^t_{tj})\right ), \label{ppMultiSchwar1}
\end{array}
\end{equation}
where we have put
$$
{}^\natural\! B^k_{ij}:=\left (A^{kr}_{j}\,\omega_i-A^{r}_{ij}\,
\omega^k -A^{k}_{sj}\, A^{sr}_i+A^{rk}_u\, A^u_{ij}\right ) \, d_r
(\log F),
$$
to avoid clutter; ${}^\natural\! D$ is the covariant derivative
associated with the Cartan connection, $A_{ijk}$ are the
components of the Cartan tensor (\ref{car}), $\omega_i$ are the
components of the Hilbert form (\ref{hil}),
${}^\natural\!\ell(f)_{ij}^k$ are the components of the tensor
above, $d_r:=u\circ d$ and $u$ is the map as in (\ref{pi}).

Theorem \ref{main} still holds for the operator ${}^\natural\!
{\cal A}(f).$
\begin{rmk}{\rm (i) The operator ${}^\natural \!{\cal A}(f)$ written
by means of the Cartan connection coincides with the operator
${\cal A}(f)$ written by means of the Chern connection only and
only when the Finsler function $F$ is Riemannian. Indeed, the
Christoffel symbols of the Cartan connection as defined in
(\ref{mus}) coincide with the Christoffel symbols of the Chern
connection only and only when the components $A^i_{jk}\equiv 0.$

(ii) the operator $\cal A$ can be expressed in terms of the
Hashiguchi connection as well. We omit here its explicit
expression. Its worth noticing that the operator $\cal A$ written
by means of the Cartan connection coincides with the operator
$\cal A$ written by means of the Hashiguchi connection when the
Finsler manifold is a Landsberg space. }
\end{rmk}
\section{Conformally invariant quantization by means of a Sazaki type metric}
\subsection{Conformally invariant quantization}
Let $(N,\am)$ be a Riemannian manifold of dimension $m.$ Denote by $\nabla$
the Levi-Civita connection associated with the metric $\am.$ We recall the
following theorem.
\begin{thm}[\cite{do}] For $m>2$ and for all $\d :=\mu-\l \not\in
\{\frac{2}{m},\frac{m+2}{2m}, \frac{m+1}{m},\frac{m+2}{m}\},$ there exists an
isomorphism
$$
{\cal Q}_{\l,\mu}^{\am}: \Pol_\d^2 (T^*N) \rightarrow \cDp
^2(\cF_{\l}(N), \cF_{\mu}(N)),
$$ given as follows: for all~$P \in~\Pol_\d^2 (T^*N),$ one can
associate a linear differential operator given by
\begin{eqnarray}
{\cal Q}_{\l,\mu}^{\am} (P)&=&
P^{ij}\, \n_i\, \n_j \nonumber\\[2mm]
&&+(\beta_1 \, \n_i\, P^{ij}+\beta_2\,\am^{ij}\, \n_i\,(\am_{kl}
P^{kl}))
\n_j
\label{Tensor1}\\[2mm]
&&+\beta_3\, \n_i\, \n_j\, P^{ij}+\beta_4\, \am^{st}\,\n_s\,\n_t\,
(\am_{ij}P^{ij})
+\beta_5\, R_{ij}P^{ij}+\beta_6 \,R \,\am_{ij}\,P^{ij},
\nonumber
\end{eqnarray}
where $P^{ij}$ are the components of $P$ and $R_{ij}$
(resp. $R$) are the Ricci tensor components (resp. the scalar curvature) of
the metric $\am;$ constants $\beta_1,\ldots,\beta_6$ are given by
\begin{equation}
\label{con}
\begin{array}{l}
\beta_1=\displaystyle
\frac{2(m\l+1)}
{2+m(1-\d)},\\[3mm]
\beta_2= \displaystyle
\frac{m(2\l+\d-1)}
{(2+m(1-\delta ))
(2-m\d)},\\[3mm]
\beta_3= \displaystyle
\frac{m\l(m\lambda+ 1)}
{(1+m(1-\delta))(2+m(1-\delta ))},\\[3mm]
\beta_4= \displaystyle
\frac{m\l(m^2 \mu (2-2\l-\d )+2(m\l +1)^2-m(m+1))}
{(1+m(1-\delta))(2+m(1-\d))(2+m(1-2\d))(2-m\delta )}, \\[3mm]
\beta_5= \displaystyle
\frac{m^2\l(\l+ \d -1)}
{(m-2)(1+m(1-\delta))} ,\\[3mm]
\beta_6= \displaystyle
\frac{(m\d-2)}
{(m-1)(2+m(1-2\delta))}\, \beta_5\cdot
\end{array}
\end{equation}
The quantization map ${\cal Q}^{\am}_{\l,\mu}$ has the following properties:\\
(i) it does not depend on the rescaling of the metric $\am$;\\
(ii) if $N=\bbR^m$ and it is endowed with a flat conformal structure, this map
is unique, equivariant with respect to the action of the group
$\Og (p+1,q+1) \subset \Diff(\bbR^m),$ where $p+q=m.$
\end{thm}
\subsection{A Sazaki type metric on $TM\backslash 0$ }
Let $(M,F)$ be a Finsler manifold of dimension $n.$
The Finsler function $F$ gives rise to a Sazaki (type) metric $\w$ on the
manifold $TM\backslash 0,$ given in local coordinates $(x^i,y^i)$ by
\begin{equation}
\label{met} \w:=\left
(\g_{ij}+\g_{st}\frac{N^s_i}{F}\frac{N^t_j}{F}\right ) dx^i\otimes
dx^j +\g_{is}\frac{N^s_j}{F^2} dy^i\otimes
dx^j+\g_{is}\frac{N^s_j}{F^2} dx^j\otimes dy^i+
\frac{\g_{ij}}{F^2}d y^i\otimes d y^j,
\end{equation}
where $\g_{ij}$ are the components of the fundamental tensor and $N^i_j$ are
given as in (\ref{en}).
\begin{rmk}
{\rm
Let us emphasize the difference between the geometric objects $\w$ and $\g.$
The metric $\w$ is a Riemannian metric on the bundle $TM\backslash 0,$ whereas
$\g$ defines a section of the bundle $\pi^*(T^{*}M)\otimes \pi^*(T^{*}M).$
When (and only when) $F$ is Riemannian, $\g$ has no $y$ dependence and
can then descend to a metric on the manifold $M.$
}
\end{rmk}
\begin{lem}
\label{lemm1}
Any tensor $P$ on $S\Gamma (TM^{\otimes 2}) $ can be extended to
a tensor on $S\Gamma(T(TM\backslash 0) ^{\otimes 2}),$ given in local
coordinates $(x^i,y^i)$ by
\begin{equation}
\label{av}
\nonumber
\tilde P=P^{ij}\frac{\partial}{\partial x^i}\frac{\partial}{\partial x^j}
- P^{it}N^j_t\frac{\partial}{\partial x^i}\frac{\partial}{\partial y^j}
- P^{js}N^i_s\frac{\partial}{\partial y^i}\frac{\partial}{\partial x^j}
+\left (P^{st}N_s^iN_t^j+P^{ij}F^2\right )
\frac{\partial}{\partial y^i}\frac{\partial}{\partial y^j},
\end{equation}
where $P^{ij}$ are the components of the tensor $P$ and $N^i_j$ are given as
in (\ref{en}).
\end{lem}
{\bf Proof.}
The objects $N^i_j,$ $\partial / \partial x^i $ and $\partial/ \partial y^i,$
behave under coordinates change as follows: for local changes on $M$,
say $(x^i)$ and their inverses $(\tilde x^i),$ one has
\begin{equation}
\label{tra}
\begin{array}{lcl}
\tilde N^i_j&=&\displaystyle  \frac{\partial x^t}{\partial \tilde x^j}
\frac{\partial \tilde x^i}{\partial  x^s}N^s_t+
\frac{\partial \tilde x^i}{\partial x^t}
\frac{\partial^2 x^t}{\partial \tilde x^s \partial \tilde x^j}
\tilde y^s,\\[3mm]
\displaystyle \frac{\partial }{\partial \tilde x^i}&=&\displaystyle
\frac{\partial x^p}{\partial \tilde x^i}
\frac{\partial}{\partial x^p}
+\frac{\partial^2 x^j}{\partial \tilde x^i \partial \tilde x^s}\tilde y^s
\frac{\partial}{\partial y^j},\\[3mm]
\displaystyle  \frac{\partial}{\partial \tilde y^i}&=&
\displaystyle \frac{\partial  x^p}
{\partial \tilde x^i}\frac{\partial}{\partial y^p}\cdot
\end{array}
\end{equation}
By substituting these formulas into (\ref{av}) we see that the
geometrical object $\tilde P$ behaves under coordinates change as
a symmetric twice-contravariant tensor field on $TM\backslash 0.$
\cqf
\begin{lem}
\label{lemm2}
The space $\cF_{\l}(M)$ can be identified with the subspace of
$\cF_{\frac{\l}{2}}(TM\backslash 0)$ with elements of the form
\begin{equation}
\label{ara}
\nonumber
\phi(x)
\left (dx^1\wedge \cdots \wedge dx^n\wedge d y^1 \wedge
\cdots \wedge d y^n
\right )^{\frac{\l}{2}}.
\end{equation}
\end{lem}
{\bf Proof.}  The 1-forms $dx^i$ and $dy^i$ behave under coordinates change
in $TM\backslash 0$ as follows: for local coordinates change on $M$,
say $(x^i)$ and their inverses $(\tilde x^i),$ one has
$$
\begin{array}{lcl}
d \tilde y^i&=&\displaystyle
\frac{\partial \tilde x^i}{\partial x^p} d y^p
+\frac{\partial^2 \tilde x^i}{\partial x^s \partial  x^t}y^s
dx^t,\\[3mm]
d \tilde x^i&=&
\displaystyle \frac{\partial  \tilde x^i}
{\partial  x^p}d x^p\cdot
\end{array}
$$
By substituting these formulas into (\ref{ara}) we see that the
geometrical object (\ref{ara}) behaves under coordinates change as
a tensor density of degree $\l/2$ on $TM\backslash 0.$
\cqf
\\

The second main result of this paper is to give some properties of the
conformally invariant quantization map (\ref{Tensor1}) by means of the Sazaki
(type) metric (\ref{met}). Namely
\begin{thm} For any lift of $P\in \Pol_{2\d}^2 (T^*
M),$ the quantization map
\begin{equation}
\label{met2} {\cal Q}^{\w}_{\l,\mu}: \Pol_\d^2(T^*(TM\backslash
0))\rightarrow \cDp ^2(\cF_{\l}(TM\backslash
0),\cF_{\mu}(TM\backslash  0)),
\end{equation}
has the property that the operator ${\cal Q}^{\w}_{\l,\mu}(\tilde
P)$ cannot descend as an operator acting on the space of
differential operators on tensor densities on $M.$ However, when
the Finsler function $F$ is Riemannian, viz
$F=\sqrt{\g_{ij}y^iy^j},$ three properties are distinguished:

(i) If $\l\not =0$ and $\mu\not = 1,$ the operator
${\cal Q}^{\w}_{\l,\mu}(\tilde P)$ cannot descend;

(ii) If $\l=0$ and $\mu\not=1$ (or $\l\not=0$ and $\mu=1$), the
operator ${\cal Q}^{\w}_{\l,\mu}(\tilde P)$ can descend only if
$\g$ is the Euclidean metric;

(iii) if $(\l,\mu)=(0,1),$ the operator ${\cal
Q}^{\w}_{0,1}(\tilde P)$ can descend, given explicitly in terms of
the metric $\g$ by
$$
P^{ij}\,\, {}^\g\n_i\,{}^\g \n_j\,+\,{}^\g\n_i (P^{ij})\,\,
{}^\g\n_j,
$$
for all $P \in \Pol_{2\d}^2 (T^*M).$
\end{thm} {\bf Proof.}
First, Lemma (\ref{lemm1}) and (\ref{lemm2}) assure that the
operator ${\cal Q}^{\w}_{\l,\mu}(\tilde P)_{|_{\cF_{2 \l}(M)}}$ is
a well-defined operator.

Suppose that $F$ is not Riemannian. In that case, the metric
$\w$ (\ref{met}) depends, in any local coordinates $(x^i,y^i),$ on $x$ and on
$y$ as well. It is easy to see from the map (\ref{Tensor1}) that
${\cal Q}_{\l,\mu}^{\w}(\tilde P)_{|_{\cF_{2\l}(M)}}$ depends on $y.$ The
crucial point of the proof is when $(\l,\mu)=(0,1).$ In that case, let us
exhibit the operator ${\cal Q}^{\w}_{0,1} (\tilde P)_{|_{\cF_{0}(M)}}$ in
local coordinates $(x^i,y^i)$, namely
\begin{equation}
\label{kin}
{\cal Q}^{\w}_{0,1} (\tilde P)_{|_{\cF_{0}(M)}}=
P^{ij}\,\partial_{x^i}\, \partial_{x^j} +( \partial_{x^i}\, P^{ij}
-P^{sj}\, \partial_{y^i}\,(N^i_s )) \, \partial_{x^j},
\end{equation}
where the components $N_s^i$ are given as in (\ref{en}). Since $F$ is not
Riemannian, the components $\partial_{y^i}( N^i_s )$ still have $y$ dependance.

Suppose now that $F$ is Riemannian. Let us prove (i) and (ii)
simultaneously. Suppose, without lost of generality, that
$M=\bbR^n$ and $\g$ is the Euclidean metric. Namely,
$\g:=\d_{ij}\, dx^i\otimes dx^j,$ where $(x^i)$ are local
coordinates on $\bbR^n.$ To achieve the proof, we will express the
quantization map (\ref{met2}) in these local coordinates, and
prove that it has $y$ dependance.

In the coordinates mentioned above, the Christoffel symbols of the
Levi-Civita connection associated with the metric $\w$ (\ref{met}) are given by
\begin{equation}
\label{lift}
\begin{array}{ll}
{\Gamma_{ij}^k}^{\w} =0
&{{\Gamma}_{i{\bar j}}^{k}}^{\w}= 0,\\[3mm]
{{\Gamma}_{{\bar i}{\bar j}}^{ k}}^{\w}=0,&
{\Gamma_{ij}^{\bar k}}^{\w}=0,
\\[3mm]
{\Gamma_{i\bar j}^{\bar k}}^{\w} =0, & \displaystyle
{\Gamma_{{\bar i}\bar j}^{\bar k}}^{\w}=-\frac{1}{F}\left (
\omega_j\,\delta^k_i
+\omega_i\,\delta^k_j
-\omega^k\,\delta_{ij}\right),
\end{array}
\end{equation}
where $\omega_i$ are the components of the Hilbert form (\ref{hil}). \\
It follows that the contraction between the tensor
$\tilde P\in \Pol_\d (T^*(TM\backslash 0))$ with the Ricci tensor $Ric$ of
the metric $\w$ is given by the equations
\begin{equation}
\label{ric}
\begin{array}{ll}
R_{ij}\, P^{ij}=0, &R_{\bar ij}\, P^{\bar i j}=0, \\[2mm]
R_{i\bar j}\, P^{i\bar j}=0, &
R_{\bar i\bar j}\, P^{\bar i\bar j}=(n-2)\,\omega_i\,\omega_j\, P^{ij}
+(2-n)\,\g_{ij}\, P^{ij},
\end{array}
\end{equation}
where $R_{ij}$ are the components of the Ricci tensor $Ric.$ \\
A direct computation proves that the scalar curvature of the metric $\w$ is
equal to
\begin{equation}
\label{scal}
3n-n^2-2.
\end{equation}
Now we are in position to express the quantization map (\ref{met2}) in local
coordinates. Using the explicit formul{\ae} of the connection (\ref{lift}),
the formul{\ae} (\ref{ric}) and (\ref{scal}), we will see that the
quantization map ${\cal Q}_{\l,\mu}(\tilde P)$ restricted $\cF_{2\l}(M)$ turns
out to be of the form
$$
\begin{array}{l}
P^{ij}\, \partial_{x^i} \, \partial_{x^j}  \\[2mm]
+\left (\, \beta_1\, \partial_{x^i}\, P^{ij} +
2\, \beta_2 \, \g^{ji}\, \g_{st}\,\partial_{x^i}
P^{st}\, \right )\, \partial_{x^j} \\[2mm]
+\beta_3\, \partial_{x^i}\, \partial_{x^j} \, P^{ij}+
2\, \beta_4 \,\g_{ij}\, \partial_{x^s} \, \partial_{x^t}\,P^{ij}
+ \left(\, (3n-n^2-2)\, \beta_6 +(2-n)\, \beta_5\,\right)
\,\g_{ij}\, P^{ij}\\[2mm]
+\displaystyle 3\,\frac{n^2 \l (\mu-1)}{1+2n}\, \omega_i
\,\omega_j\, P^{ij},
\end{array}
$$
where constants $\beta_1,\ldots, \beta_6$ are given as in
(\ref{con}). As the constant $\frac{n^2 \l (\mu-1)}{1+2n}$ does
not vanish when $\l\not=0$ and $\mu\not=1$, the quantization map
still have $y$ dependence, and then does not take it is values in
$\cF_{2\mu}(M).$ Part (i) is proven. To achieve the proof of part
(ii), let us consider a non-Euclidean metric $\g.$ In local
coordinates $(x^i,y^i)$, the operator ${\cal
Q}_{\l,\mu}^{\w}(\tilde P)$ restricted to $\cF_{2\l}(M)$ will have
a component of the form
$$(1-\beta_1)\left (\frac{1}{F^2}\,\partial_{x^j} (N_s^u)\, N_i^v\,
\g_{uv}\,\g^{ks}\,P^{ij}\right )\partial_{x^k},
$$ which has $y$ dependance. The constant
$1-\beta_1$ does not vanish under the condition $\l=0$ and $\mu\not =1$
(or $\l\not=0$ and $\mu =1$). Part (ii) is proven.

Let us prove part (iii). For $(\l,\mu)=(0,1),$ the quantization map, written
in any local coordinates $(x^i,y^i)$, has the form
\begin{equation}
\label{ekl}
P^{ij}\,\partial_{x^i} \partial_{x^j} +\left ( \partial_{x^i} P^{ij}
-\frac{1}{2}\, \g^{uv}\, \partial_{x^i}\, \g_{uv}\,P^{ij} \right )
\partial_{x^j}\cdot
\end{equation}
Using the following formul{\ae}:
$$
\begin{array}{lcl}
P^{ij}\,\partial_{x^i}\,\partial_{x^j}&=&P^{ij}\,\,
{}^\g\n_i{}^\g\, \n_j
+ {}^{\g}\Gamma^i_{jk}\, P^{jk}\,\,{}^\g\n_i, \\[3mm]
\partial_{x^i}\,P^{ij}\,\partial_{x^j}&=&{}^\g \nabla_i P^{ij}\,\,{}^\g \n_j
 + {}^{\g}\Gamma^i_{jk}\, P^{jk}\,\, {}^\g\n_i
 +P^{ij}\,\,{}^\g\Gamma_j\,\,{}^\g \n_i,
\end{array}
$$
we see that the formula (\ref{ekl}) turns into the form
$$
P^{ij}\,\,
{}^\g\n_i{}^\g\, \n_j+{}^\g \nabla_i P^{ij}\,\,{}^\g \n_j. \\
$$
The formula above has certainly no $y$ dependance. Part (iii) is
proven. \cqfd
\begin{rmk}
{\rm Theorem above shows that the quantization map ${\cal
Q}^{\w}_{\l,\mu}(\tilde P)_{|_{\cF_{2\l}(M)}}$ does not reproduce
the quantization map ${\cal Q}^{\g}_{2\l,2\mu}(P)$ even if $F$ is
Riemannian. }
\end{rmk}
\section{Open problems}
{\bf 1)} Following \cite{b1}, there exists two 1-cocycles on the
group $\Diff(M)$, say $c_1,c_2,$ with values in $\cDp (S\Gamma
(TM^{\otimes 2})\otimes \cF_{\d}(M), \Gamma (TM)\otimes
\cF_{\d}(M))$ and $\cDp (S\Gamma (TM^{\otimes 2})\otimes
\cF_{\d}(M),\cF_{\d}(M))$ respectively, that are conformally
invariant; namely, they depend only on the conformal class of the
Riemannian metric. These 1-cocycles were introduced in \cite{b1}
as conformal multi-dimensional Schwarzian derivatives. In this
paper, we have introduce the 1-cocycle $\cal A$ (see
(\ref{MultiSchwar1})) as the Finslerian analogous to the 1-cocycle
$c_1;$ however, the computation to extend the 1-cocycle $c_2$
seems to be more intricate.\\

{\bf 2)} We ask the following question:\\
{\it Is there a map
$$
{\cal Q}: \Pol (T^*(TM\backslash 0))\otimes
\cF_{\mu-\lambda}(\pi^*(T^*M)) \rightarrow {\cal D}
(\cF_\l(\pi^*(T^*M)),  \cF_\mu(\pi^*(T^*M))),
$$
having the following properties:

(i) it does not depend on the rescaling of the Finsler function by
a non-zero positive function on $M;$

(ii) it coincides with the Duval-Ovsienko's conformally invariant
map when $F$ is
Riemannian?}\\

We believe that a positive answer to this question will probably
produce the 1-cocycle $c_2$ discussed in part {\bf 1)}. It should
be stressed, however, that the quantization map and the 1-cocycle
$c_2$ may not exist in the {\it generic} Finsler setting. Such a
situation happened in Conformal Geometry where a large number of
invariant differential operators do not generalize to arbitrarily
``curved'' manifolds. For example, the power of the Laplacian,
$\Delta^k,$ where $k$ is an integer, are the unique differential
operators acting on the space of tensor densities of appropriate
weights, that are invariant under the action of the conformal
group (see \cite{er, jv}); however, their curved analogues do not
exist when the dimension of the manifold $\rm{dim }\, M$ is
greater than $4$ and even, and $k> \rm{dim }\, M/2$, as recently
proven in \cite{gh}.

\bigskip

{\it Acknowledgments}. It is a pleasure to acknowledge numerous
fruitful discussions with A. Cardona, Ch. Duval, V. Ovsienko and
M. Stienon.

\vskip 1cm


\end{document}